\documentclass[12pt]{article}
\usepackage{amssymb}

\begin{document}
\def\square#1{\vbox{\hrule
\hbox{\vrule\hbox to #1 pt{\hfill}\vbox{\vskip #1 pt}\vrule}\hrule}}
\newtheorem{guess}{Proposition }[section]
\newtheorem {theorem}[guess]{Theorem}
\newtheorem{lemma}[guess]{Lemma}
\newtheorem{corollary}[guess]{Corollary}
\newtheorem{example}[guess]{Example}
\newtheorem{remark}[guess]{Remark}
\newtheorem{definition}[guess]{Definition}

\centerline{\sc \large  Vector bundles on a three dimensional  } 
\vspace{3mm}
\centerline{\sc \large  neighborhood of a ruled surface }
\vspace{7mm}
\centerline{\large Edoardo Ballico \hspace{1.5 cm} Elizabeth Gasparim} 
\centerline{ \,\,\,\, \, University of Trento \hspace{1.6cm} New Mexico State University}      
\begin{abstract} Let $S$ be a ruled surface inside a smooth 
threefold $W$
and let $E$ be a vector bundle on a formal neighborhood of $S.$ 
We find minimal conditions under which the local moduli space 
of $E$ is finite dimensional and smooth. Moreover, we show that 
$E$ is a flat limit of a flat family of vector bundles whose 
general element we describe explicitly.
\end{abstract}

\section{Introduction}
Consider the general question: how do moduli spaces of vector bundles
change under birational transformations of the base?
In this paper, we take the first steps of a program to 
study this question for threefolds.
In dimension three, flops give essential examples of 
birational transformations. 

We first recall the definition of the basic flop.
Let $X$ be the cone over the ordinary double point defined by the equation 
$xy-zw=0$ on ${\mathbb C}^4.$
The basic flop is described by the diagram:
$$  \begin{matrix}{ &  \widetilde{X}  & \cr
                   {\,}^{f_1}\!\! \swarrow
                   & &  {\searrow}^{\!f_2}  \cr
                   X_1& & X_2  \cr
                   {\,}_{\pi_1} \!\! \searrow 
                   & & {\swarrow}_{\!\pi_2}\cr
                  & X &} \end{matrix} $$
where 
$\widetilde{X}\colon = \widetilde{X}_{x,y,z,w}$ is 
the blow up of X at the vertex $x=y=z=w=0,$
$X_1\colon =\widetilde{X}_{x,z}$ is the blow up of X along  $x=z=0$
and $X_2\colon = \widetilde{X}_{y,w}$ is the blow up of X along  $y=w=0.$
The  {\em basic flop}  is the transformation 
from $X_1$  to $X_2.$ 
The spaces appearing in this diagram are not compact, but they 
do contain neighborhoods of compact curves. 
We wish to find what vector bundles fit over this
 diagram, together with  their local deformations.
Note that on the given diagram, the spaces $X_1$ and $X_2$ are both 
abstractly 
isomorphic to ${\cal O}_{{\mathbb P}^1}(-1)  \oplus 
{\cal O}_{{\mathbb P}^1}(-1),$
although the maps $\pi_1 $ and $\pi_2$ are distinct; whereas
$\widetilde{X}$ is isomorphic to   ${\cal O}_{{\mathbb P}^1 \times
{\mathbb P}^1}(-1,-1),$

 We generalize the situation  of $\widetilde{X}$ slightly by
considering a  ruled surface $S$ 
with negative normal bundle
inside a smooth threefold.
We  then  study bundles $E$ on a formal neighborhood 
$\widehat S$ of 
$S$ and   their local moduli spaces (cf. definition \ref{local}).
In the case of a Hirzebruch surface $S,$ 
we require that $E\vert_S$ be simple. 
When $S$ is ruled over a curve $C$ of genus greater than 1,
we assume that $E\vert_S$ is $R$--stable with respect to  a
good polarization $R$ of $S$ (cf. definition \ref{def}).

These conditions are minimal in the following sense. The local moduli 
space of a simple bundle is unobstructed, and therefore smooth (cf. remark 
\ref{smooth}). 
Hence to have smoothness it would be desirable to impose the condition that 
$E$ be simple. 
However, as an easy argument in section 3 shows,
 there are no simple bundles on 
$\widehat S.$  The alternative is to impose a condition 
on the restriction of $E$ to an infinitesimal neighborhood
of $S.$ We choose the zero-th formal neighborhood. 
We have the following results.
\vspace{5mm}

\noindent {\bf Theorem A} {\it  Let 
$S$ be a Hirzebruch surface
with negative normal bundle inside a smooth threefold.
   Let $E = \{E_n\}$  be a vector bundle on $\widehat{S}$ 
such that $E\vert_S$ is simple. Then the local moduli space of $E$ is
 finite dimensional and smooth. Moreover, $E$ is a flat limit of a flat
 family 
of vector bundles on $\widehat{S}$  satisfying  
 properties $(\iota)$ and $(\iota\iota)$ below.}

\vspace{5mm}
\noindent {\bf Theorem B}
 {\it Let $S$ be a ruled surface with negative normal bundle inside a 
smooth threefold, so that $S$ is ruled over a curve of positive genus.
   Fix a good polarization $R$ on $S.$
 Let $E =\{E_n\}$  be a vector bundle on
$\widehat{S}$ such that $E\vert_S$ is $R$--stable. 
Then the local moduli space of $E$ is finite 
dimensional and smooth. Moreover, $E$ is a flat limit 
of a flat family of vector bundles on $\widehat{S}$ 
satisfying  properties $(\iota)$ and $(\iota\iota)$ below.}

\vspace{5mm}
\noindent  Let $r\colon=
\mbox{rank}(E)$ and $d\colon = \mbox{deg}(E).$
The general element $ G= \{G_n\}$ of the family has the following behavior.

\vspace{3mm}

\noindent{\it ($\iota$)\, 
If $d=ar-x, \,0 < x < r, $ then the general 
element $G$ of the family   is a vector bundle such that
the restriction of $G_1$ to a general fiber $ D$ 
of $u$ has  splitting type $(a,\cdots ,a,a-1,\cdots ,a-1),$ 
and in this case $$G|_{\widehat{D}} \simeq
 {\cal O}_{\widehat{D}}(a)^{\oplus(r-x)}
\oplus {\cal O}_{\widehat{D}}(a-1)^{\oplus x }.$$}

{\it \noindent($\iota\iota$)\, 
If $d=ra,$ then the general element $G$ of the family 
is a vector bundle such that the restriction 
 of $G_1$ to a general fiber $D$ of $u$ has splitting type $(a, \cdots ,a)$
and in this case  $$G|_{\widehat{D}}
\simeq {\cal O}_{\widehat{D}}(a)^{\oplus r}$$
but there exists a finite number of jumping fibers $D'$ 
where $G|_{D'}$ has splitting type  $(a+1,a,...,a,a-1).$
 }
\vspace{3mm}

For a bundle over a Hirzebruch surface
 we calculate the number of such jumping fibers.

\vspace{5mm}
\noindent {\bf Theorem C} {\it 
Let $z$ be number of  jumping fibers of $G.$
Set  $E=G(-a{\bf h})$ and $m\colon=  \mbox{deg}(u_*E).$ Then 
$$z=c_2 (E)=c_2(G)-a(r-1)c_1(G)\cdot {\bf h}-ea^2r(r-1)/2$$
 and 
$$m=c_1(u_*E)=-z+c_1(G)\cdot{\bf h}+rae.$$}

\noindent In section 2 we recall some basic concepts of deformation 
theory.
In section 3 we consider bundles 
on a neighborhood of a Hirzebruch surface and 
prove Theorems A and C. In section 4 we consider bundles on a 
neighborhood of a surface ruled over a curve of higher genus 
and prove Theorem B.

\section{Background material on deformations}

In this paper we work only over $\mathbb C.$
The basic material on the deformation theory appearing in this section is
taken from Seshadri \cite{Se}.  
Let $X$ be a
 scheme over and algebraically
closed field $k.$
Let $R$ be a complete local ring such that $R/m_R=k$, 
$m_R$ the maximal ideal of $R$ and $R_n= R/m_R^n.$

\begin{definition}{\em 
A {\it  deformation } $Y$ of $X$ parametrized by a scheme $T$ 
with base point $t_0$ consists of 
\begin{enumerate}
\item a morphism $Y \rightarrow T$ which is flat and of
finite type
\item a closed point $t_0 \in T,$ and an isomorphism 
$Y_{t_0}  \stackrel{\sim } {\rightarrow }X,$
where  $Y_{t_0}= Y\times _T k(t_0)$
is the fiber over $t_0.$ 
\end{enumerate}}
\end{definition}

\begin{definition}{\em
 A {\it formal deformation }
$X_R$ of $X$ is a sequence $\{X_n\}$ such that
\begin{enumerate}
\item $X_n=X_{R_n}$ where $X_{R_n}$  is a deformation of $X$ over 
$R_n$ 
\item we are given a compatible sequence of isomorphisms 
$X_n \otimes _{R_n} R_{n-1} \rightarrow X_{n-1}$ for any $n.$
\end{enumerate}}
\end{definition}

\begin{definition}{\em
 Let $A$ be a finite dimensional local $k$-algebra. 
Then, giving a $k$-algebra homomorphism 
$\phi \colon R \rightarrow A$ is equivalent to  giving a 
compatible sequence of homomorphisms $\phi_n\colon R_n \rightarrow A$
 for $n >>0.$ It follows that, given a formal deformation 
$X_R$ of $X$ and a homomorphism $\phi\colon R \rightarrow A,$ 
$X_n \otimes_{R_n} A$ is the same up to isomorphisms for $n>>0.$ We 
define this to be $X_R \otimes A.$ 
 It is a deformation of $X$ over $A,$ called the 
{\it base change} of $X_R $ by $\mbox{Spec } A  \rightarrow  \mbox{Spec } X.$ 
}
\end{definition}

\begin{definition}{\em 
Let $F$ and $G$ be the functors
 defined by
\begin{enumerate}
\item $F(A)$ = isomorphism classes of deformations $X_A$ over $A$
\item $G(A)= Hom _k(R,A).$
\end{enumerate} We get a morphism of functors 
$j\colon G \rightarrow F$ defined by 
$\phi \in  Hom _k (R,A) \mapsto X_R \otimes A.$
A formal deformation $X_R$ of $X$ is said to be {\it versal} if the 
functor $j$ is
formally smooth. (\cite{Se} p.271)}
\end{definition}
More generally, one can define similarly the concept of versal deformation
for a covariant functor $F$ with $F(k)$ = a single point. 
Schlessinger gave conditions for the existence of  versal 
deformations of a functor for $F.$ Moreover, Artin's 
algebraization theorem says that Schlessinger's conditions 
together with effectiveness imply the 
existence of an  algebraic deformation space for  $F.$
For details see \cite{Ar} and \cite{Sc} .

\begin{remark}{\em
In the case of deformations of $X$  
algebraization  means that there exists a 
scheme $Y$ over $S$ flat and of finite type, with base point $s_0$ such that 
$\widehat{{\cal O}}_{S,s_0}=R$ and $Y \otimes R_n = X_n.$ 
The conditions for algebraization are satisfied 
for deformations of vector bundles over a complete algebraic scheme 
(see \cite{Se} thm 2.3). }
\end{remark}

\begin{definition}\label{local}{\em
The germ  
of $Y$ at $s_0 \in S$ is determined 
up to isomorphism and we call it the {\it local moduli space}
of $X.$ (Here germ  means  Spf ${{\cal O}}_{Y, s_0}  ,$ i.e. a 
compatible sequence of spectra of rings over Artinian rings).
When a deformation of $X$ is considered only on
a germ at a point $s_0$ we call it a {\it local deformation}
of $X.$ 

% (in the case of formal deformations of a functor 
%$F$ Schlessinger called it  the hull of $F.$)  
}\end{definition}
 
In this paper we construct formal
deformations of vector bundles. An application of 
Artin's algebraization theorem then
implies the existence of algebraic
deformations. We restate definition 2.1  for 
the case of vector bundles (definitions 2.2 and 2.4 can be repeated 
with $E$ in place of $X$ and $\overline{E} $ in place of $Y$ and we get the 
concepts of formal, versal and local deformations as well as local moduli
for vector bundles).

\begin{definition}{\em  Let $E$
be a vector bundle over a scheme $X.$
A {\it  deformation } $\overline{E}$ of $E$ parametrized by a scheme $T$ 
with base point $t_0$ consists of 
\begin{enumerate}
\item a vector bundle  $\overline{E} \rightarrow X \times T$
\item a closed point $t_0 \in T,$ and an isomorphism 
$\overline{E}\vert_{X \times t_0}  \stackrel{\sim } {\rightarrow }E.$
\end{enumerate}}
\end{definition}

\begin{definition}{\em Let $X_R = \{X_n\}$ be a formal deformation of 
$X.$ A {\it vector bundle} $E_R$ on $X_R$ is a compatible sequence of 
vector bundles on each $X_n.$
 A {\it  deformation } of 
$E_R$ is given by a compatible sequence of deformations 
for each $E_n.$ }
\end{definition}

\begin{remark}\label{smooth}{\em
We say that the local moduli space of $E$ is {\it smooth} if the germ 
of $\overline{E}$ at $X\times t_0$ is regular.
In order to check that the local moduli space of $E$ is smooth  it 
suffices to check formal smoothness (cf. \cite{Se} remark 2.4).
Obstructions for smoothness are in $H^2(X,End(E)),$ 
and if this group vanishes, we say that  deformations of $E$
are {\it unobstructed}. 
It follows that the criterion for smoothness 
is that $H^2(X,End(E)) = 0.$  }
\end {remark}

\section{Bundles on  neighborhood of a Hirzebruch surface} 
 
Let $S$ be a ruled surface inside a smooth threefold $W.$
Let  $V$ be either a neighborhood of  $S$ in $W$
in the smooth topology, or the germ of $W$ around $S,$  and let 
$\widehat{S}$ be the formal completion of $S$ in $V.$
In this section we consider the case when 
 $S = \Sigma_e, e\geq 0,$ is a Hirzebruch surface.
 If $e = 0$ then $S \simeq {\mathbb P}^1 \times {\mathbb P}^1$
 and hence the two projections $f_1\colon S \rightarrow {\mathbb P}^1$ 
 (resp. $f_2\colon S \rightarrow{\mathbb P}^1 )$ on the first (resp. second)
 factor define two rulings of $S.$ 
We use only the first ruling and set $u\colon = f_1.$
 If $e > 0,$ then the surface $\Sigma_e$ has a unique ruling
 $u\colon \Sigma_e \rightarrow {\mathbb P}^1 .$
  Call $ f$ any fiber of $u.$ 
Let $h$ be a section of $u$ with minimal self-intersection.
 We  denote by $\bf f$ (resp. $\bf h)$
 the class of $f$
 (resp. $h)$ in Pic$(S).$ Thus $\bf f$ and $\bf h$
 form a basis for  Pic$(S) \simeq {\mathbb Z} \oplus {\mathbb Z}$
and have
 intersection numbers  ${\bf h}^2 = -e,\, {\bf h}\cdot {\bf f} = 1$
 and ${\bf f}^2 = 0.$
 The canonical line bundle of $S$ is isomorphic to
 ${\cal O}_S(-2{\bf h}-(e+2){\bf f}).$
Let $t$ and $s$ be the integers such that
 $I/I^2 \simeq {\cal O}_S(t{\bf h}+s{\bf f}),$ where $I$ is the ideal  
defining 
$V$ in $W.$
 We assume $t > 0$ and $s > et,$ that is,
 we assume that $I/I^2$ is ample. 
For every integer $n\geq 1$ we have
$ I^n/I^{n+1} \simeq 
S^n(I/I^2) \simeq {\cal O}_S(nt{\bf h}+ns{\bf f}).$
 Thus $h^1(S,I^n/I^{n+1}) = h^2(S,I^n/I^{n+1}) = 0$
 for all $n\geq 1.$

The vector spaces of regular  functions on $V$ 
and of  formal  functions on 
$\widehat{S}$ are infinite dimensional.
  We consider vector bundles $E = \{E_n\}$  on $\widehat{S}$
such that $E\vert_S$ is simple, that is, such that
$ h^0(S,End(E\vert_S)) = 1.$ In other words, we require
 $ h^0(S,ad(E\vert_S)) = 0.$
 For all integers $n\geq 0$ we have the exact sequence
$$0 \rightarrow I^n/I^{n+1} \rightarrow {\cal O}_{S(n+1)} 
\rightarrow {\cal O}_{S(n)} \rightarrow 0	\eqno{(1)}	$$
For every vector bundle $E = \{E_n\}$  on $\widehat{S}$ 
we have the exact sequences
$$ 0 \rightarrow E_{0} \otimes I^n/I^{n+1} \rightarrow E_{n+1}
 \rightarrow E_n \rightarrow 0		\eqno{(2)}	$$
obtained from $(1)$ by tensoring with $E_{n+1}.$
 Take a vector bundle $G = \{G_n\}$  on $ \widehat{S}.$
 For every integer $n\geq 0,$ set $E_n = ad(G_n),$
 where the $ad$ and  $Hom$ functors are computed with respect
 to ${\cal O}_{S(n)}.$ By the long exact sequence
in cohomology derived from  $(2)$
 we see that the integer $h^0(S(n),ad(G_n))$
 goes to infinity when n goes to infinity. Hence, there are no 
simple vector bundles on $\widehat{S}.$

\begin{lemma} \label{1} Let $E = \{E_n\}$
  be a vector bundle on $\widehat{S}$
such that $E\vert_S$ is simple. Then for all integers 
$n\geq 1$ we have $h^2(S(n),End(E_n)) = 0.$
\end{lemma}

\noindent{\em Proof}.
First assume $n = 1.$ Since $E\vert_S$ is simple, we have
$ h^0(S,End(E\vert_S)\otimes{\cal O}_S(-2{\bf h}-(e+2){\bf f})) = 0.$
 By Serre duality  
$$h^0(S,End(E\vert_S)\otimes{\cal O}_S(-2{\bf h}-(e+2){\bf f})) =
 h^2(S,End(E\vert_S)),$$ concluding the case $n = 1.$ 
Now assume $n\geq 2$ and that the result is true for the integer $n-1,$ 
i.e. assume $h^2(S(n-1),End(E_{n-1})) = 0.$ 
Since dim$S(n) = 2$ we have $h^3(S(n),A) = 0$ 
for every coherent analytic sheaf $A$ on $S(n).$ 
Using $(2)$ for the integer $n-1$ and the vector bundle
 $End(E_n)$ together with the inductive assumption, we see 
that $ h^2(S(n),End(E_n)) = 0$ if 
$h^2(S,End(E\vert_S)\otimes{\cal O}_S(nt{\bf h}+ns{\bf f})) = 0.$
 By Serre duality we have $$h^2(S,End(E\vert_S)\otimes{\cal O}_S
(nt{\bf h}+ns{\bf f}))
 =$$
$$ h^0(S,End(E\vert_S)\otimes{\cal O}_S
(-(2+nt){\bf h}-(2+e+ns){\bf f})) =0.$$ \hfill \square{5}

\begin{remark}\label{2} {\em
 If $F$ is  a vector bundle on $S(n)$ such that 
$h^2(S(n),End(F)) = 0, $ then by Remark 2.8 the local moduli 
space of $F$ is smooth and has dimension $h^1(S(n),End(F)).$ 
}
\end{remark}

\begin{lemma}\label{3} Let $E = \{E_n\}$ be   a vector bundle
 on $\widehat{S}$ such that $E\vert_S$ is simple. 
Then for all integers $n\geq 1$ the restriction map 
$$h^1(S(n+1),End(E_{n+1})) \rightarrow
 h^1(S(n),End(E_n))$$ is surjective.
\end{lemma}

\noindent{\em Proof}. As in the proof of Lemma \ref{1} we obtain
$ h^2(S,End(E\vert_S)\otimes I^n/I^{n+1}) = 0$
 for every integer $n\geq 1.$ The lemma follows from 
 the cohomology exact sequence of $(2)$ 
with the bundle $End(E_n)$ instead of $E_n.$\hfill\square{5}

\begin{lemma} \label{14} 
 Let $E =\{E_n\}$  be a vector bundle on
$\widehat{S}$ such that $E\vert_S$ is simple. 
Then there exists an integer $x$ depending only on $E\vert_S$ 
such that for all integers $n\geq x$ the restriction map 
$$H^1(S(n+1),End(E_{n+1})) \rightarrow 
H^1(S(n),End(E_n))$$ is bijective. 
\end{lemma}

\noindent{\em Proof}. By Lemma \ref{3} it  suffices to 
show the existence of $x$ such that for all $n\geq  x$ 
the restriction map $H^1(S(n+1),End(E_{n+1}))
 \rightarrow H^1(S(n),End(E_n))$ 
is injective. 
Since  $I/I^2$ is ample there exists  
 an integer $x$ 
such that for all integers $y\geq x$ we have 
$H^1(S,End(E\vert_S)\otimes I^y/I^{y+1}) = 0.$ 
Now injectivity follows from the long exact sequence 
of $(2)$ with $End(E_n)$ in place of $E_n.$
\hfill\square{5}

\begin{guess}  Let $E =\{E_n\}$  be a vector bundle on
$\widehat{S}$ such that $E\vert_S$ is simple. 
There exists an integer $x$ such that $n\geq x$ implies that 
every local deformation of $E_n$ lifts to a local deformation of 
$E.$
\end{guess}

\noindent{\em Proof}. By Schlessinger's theorem, a hull exists for 
deformations of $E_n,$
and since by lemma  \ref{1} $h^2(S(n),End( E_n))=0, $
 the hull is smooth. Hence it is the formal 
spectrum of a formal power series ring 
$R_n= \mathbb C[[x_1,\cdots, x_s]].$  Similarly, set $R_{n+1} =  
\mathbb C[[y_1,\cdots, y_r]]$ to be the formal power series ring 
corresponding to  $E_{n+1}.$ By lemma   \ref{14} the map 
$R_{n+1} \rightarrow R_n$ induces a bijection at tangent level, 
and it follows from the formal inverse function theorem, 
that the local deformations  of 
$E({n+1})$ and  $E(n)$ are isomorphic for all $n \geq x$
 therefore they determine a local deformation of $E.$ 
\hfill\square{5}
\vspace{3mm}

\noindent The following property of bundles on 
$\mathbb P^1$ is well known, but we were not able to find it in the literature.

\begin{lemma} \label{3.5}
Every  vector bundle on ${\mathbb P}^1$ 
is the flat limit of a flat family of 
rigid vector bundles.
\end{lemma}
 
\noindent{\em Proof}. Let $0 \leq c <r$ be integers and let 
$A
 = {\cal O}_{\mathbb P^1}(a)^{\oplus (r-c)}\oplus {\cal O}_{\mathbb P^1}
(a-1)^{\oplus c }.$ This is a rigid bundle. We show that all bundles 
on $\mathbb P^1$ with rank $r$ and degree $ra-c$ deform  to $A.$ 
Fix any such bundle $B$ with splitting type $b_1 \geq \cdots \geq b_r;$ 
we may assume $b_r \leq b_1-2,$ otherwise $B=A.$ 
Set $B'=  \bigoplus_{b_i \leq a-1}{\cal O}_{\mathbb P^1}(b_i)$
and  $B''=  \bigoplus_{b_i > a-1}{\cal O}_{\mathbb P^1}(b_i).$
By Shatz \cite{Sh} Proposition 1, there is an exact sequence 
$$0 \rightarrow B' \rightarrow A \rightarrow B'' \rightarrow 0.$$
Call ${\bf e}$ the extension class of this sequence. 
For all nonzero scalars $\lambda$ the extension $\lambda {\bf e}$ 
has $A$ as middle term, whereas for $\lambda =0$ 
the corresponding extension 
has middle term $B.$ This gives a flat specialization to $B$ 
of a family of vector bundles isomorphic to $A.$ 
\hfill\square{5}

\begin{lemma} \label{4.0}Let $F$ be a simple rank $r$ vector bundle on $S.$ 
 Fix a fiber $D$ of $u.$ Then the local moduli
space of the vector bundle $F|_D$ is smooth and of dimension 
$h^1(S,End(F|_D)).$ The local moduli space of 
$F$ on $S$ is smooth and of dimension 
$h^1(S,End(F)).$
\end{lemma}
 
\noindent{\em Proof}.
Since dim$(D)=1,$ we have $H^2(S,End(F|_D))=0$ and hence the local moduli 
space of the vector bundle $F|_D$ is smooth and of dimension 
$H^1(S,End(F|_D)).$
Consider 
 the exact sequence 
$$0 \rightarrow End(F)(-D) \rightarrow End(F) \rightarrow 
End(F|_D) \rightarrow 0. \eqno(3)$$
By Serre duality we have $$h^2(S, End(F)(-D))=h^0(S, End(F)(-D)\otimes 
{\cal O}_S(-2{\bf h}-(1+e){\bf f})).$$
Since $h^0(S,End(F))=1,$ we have that $h^0(End(F)(-D)\otimes 
{\cal O}_S(-2{\bf h}-(1+e){\bf f}))=0$
and hence $ h^2(End(F)(-D))= 0 .$
It follows that  $h^2(S,End(F))=0,$ we 
obtain the result for the local moduli space of $F.$ \hfill\square{5}
\vspace{5mm}

\noindent  The following observation was inspired by \cite{BH} Lemmas 2 and 3,
and their use in \cite{E}.

\begin{guess} \label{4.1} Lef $F$ be a simple vector bundle on $S.$ 
Then  $F$ is a flat 
limit of a flat family of vector bundles,  whose restriction 
to $D$ is rigid.
\end{guess}

\noindent{\em Proof}. Since $ h^2(End(F)(-D))= 0 $ as shown in the proof of 
lemma \ref{4.0},
by the exact sequence $(3)$ we obtain that the restriction map 
$$\gamma\colon H^1(S,End(F)) \rightarrow H^1(S,End(F|_D)\eqno(4)$$
 is surjective.
The surjectivity of $\gamma$ means that every local deformation of 
$F|_D$ may be lifted to a local deformation of $F.$
 Thus, we obtain from lemma \ref{3.5} that $F$ is a flat 
limit of a flat family of vector bundles,  whose restriction 
to $D$ is rigid.
\hfill\square{5}

\subsection{Case 1: the rank does not divide the degree}
 
 Let $F$ be a simple vector bundle on $S$ 
whose rank $r \geq 2$ does not divide the degree. Write
det$(F)\cdot {\cal O}_S({\bf h}) = ar-x,$ 
with $0<x<r.$
We  construct the local moduli space  
of $F.$ 

 Take any fiber $K$ of $u.$
If $F|_K$ is not rigid, then 
  by  \cite{BH} Lemma 2, it is of special type among
deformations of bundles on ${\mathbb P}^1;$ such types occurring 
in a subset of codimension at least 3 inside a complete family.
The following trick we stole from \cite{BH} and it was the starting point
of our paper.
The surjectivity of $\gamma$ in (4) implies that the general member 
of this family in Proposition \ref{4.1}
has, as restriction to any fiber, a general deformation of 
$F\vert_K,$ i.e. the only rigid vector bundle on 
${\mathbb P}^1$ with degree $ar-x, $ 
i.e. the bundle on ${\mathbb P}^1$ with
splitting type 
$(a,\cdots ,a,a-1,\cdots ,a-1)$ (with $a-1$ appearing x times).
Therefore there exists a flat deformation of $F$ 
whose general $G$ element has rigid restriction to all except finitely many 
fibers of $u.$
This implies that $G$ is uniform in the sense of Ishimura 
with respect to the ruling $u$ (see \cite{I}).
Changing bases, $H\colon =u_*(G\otimes{\cal O}_S(-a{\bf h}))$
 is a rank $r-x$ vector bundle  on ${\mathbb P}^1$  and $ u^*(H)$
 is a rank $r-x$ sub-bundle of $G\otimes{\cal O}_S(-a{\bf h})$
and  $${G\otimes{\cal O}_S(-a{\bf h})\over u*(H)} 
\otimes{\cal O}_S((a-1){\bf h}) \simeq u^*(H)$$
 for some rank $x$ vector bundle $M$ on ${\mathbb P}^1 .$ 
Thus we have an exact sequence
$$0 \rightarrow u^*(H)\otimes{\cal O}_S(a{\bf h}) 
\rightarrow G \rightarrow u^*(M)\otimes{\cal O}_S((a-1){\bf h})
 \rightarrow 0.		\eqno(5)	$$
The Chern classes of $G$ are uniquely determined
 by the integers $e, a, x, \mbox{deg}(H)$ and deg$(M).$
  Conversely, the integers deg$(M)$ and deg$(H)$ are uniquely
 determined by  $r, a, x$ and the Chern classes of $G.$
 We  set the shorthands
$$\bar M\colon = u^*(M)\otimes{\cal O}_S((a-1){\bf h}),$$
 $$\bar H\colon= u^*(H)\otimes{\cal O}_S(a{\bf h}),$$ 
$$K\colon=
 {\cal O}_S(-2{\bf h}-(e+2){\bf f}).$$
 By Serre duality  $$h^2(S,Hom(\bar M, \bar H)) = 
h^0(S, Hom (\bar M, \bar H \otimes K )).$$
 The restriction  of 
Hom$(\bar M, \bar H \otimes K)$
 to any fiber of $u$ is a direct sum of line bundles 
of degree  at most $-3.$ Thus 
$$h^2(S,\mbox{Hom}(\bar M, \bar H)) = 0. $$
By Riemann  Roch, the integer 
$h^1(S,\mbox{Hom}(\bar M,\bar H))$
 depends only on the numerical data 
 $e, r, x, a, \mbox{deg}(H),$ and deg$(M)$
 but not on the specific choice of the   vector bundles $M$ and
 $H.$ 
 From the properties  of the universal
 Ext functor (see \cite{L}) it follows  that, corresponding to  any 
family of vector bundles $\{M_t\}_{t \in T}$ and
 $\{H_t\}_{t \in T}$ parametrized by  an irreducible variety $T,$
 we obtain a flat family ${\bf V}= \{G_t\}_{t\in T}$ of
middle terms of  extensions 
$(4)$
 parametrized by $T$ with each $G_t$ a rank $r$ vector bundle
 on $T.$ 
Since T is assumed to be irreducible, ${\bf V}$ is irreducible. 

We choose  as the parameter space $T$ 
the product of a versal deformation space of $M$ 
and a versal deformation space of $H.$
 Such spaces are irreducible, smooth and of dimension respectively
$ h^1({\mathbb P}^1,End(M))$ and 
$h^1({\mathbb P}^1,End(H)).$
 With this choice of $T,$ for a general $t\in T$ 
the vector bundles $M$ and $H$ are rigid
and the the general vector bundle
$ G_t,$ of the family  is an extension of the form 
$(4)$
 with $M$ and $H$ rigid.
The set of all such extensions 
is a vector space, whose dimension depends only on the 
numerical data $r, c_1(G),$ and $c_2(G).$ 
Combining with  $(4)$
 we obtain that any simple rank $r$ vector bundle $F$ 
with degree $ar-x$ is the flat limit of a family of simple vector bundles
 $G_t$ arising from the family just described, that is,
 from an extension $(4)$ with $M$ and $H$ rigid.

We now  extend the construction of the flat family 
 to bundles on $\widehat{S}.$
Let $ D \simeq{\mathbb P}^1 $ be a fiber of $u$ and $J$ 
the ideal sheaf of $D$ in $V $ or $\widehat{S}.$ 
 Let $\widehat{D}= \{D(n)\}$
be the formal completion of $D$ in $\widehat{S}.$
% (or in $V,$ if you prefer to work in the holomorphic category,
% instead of the formal category).
Hence $D(n)$ has $J^n$ 
as ideal sheaf in $\widehat{S}$ (or in $V)$ and $D(0) = D.$ 
We have $J/J^2 \simeq {\cal O}_D(t)\otimes{\cal O}_D$ and  $J^n/J^{n+1}
 \simeq S^n(J/J^2)$ for every $n\geq 1.$ Let $A=\{A_n\}$
  be a rank $r$ vector bundle on $\widehat{S}.$
 For every integer $n\geq  1$ we have an exact sequence
$$0  \rightarrow J^n/J^{n+1}\otimes A_1 
\rightarrow A_{n+1} \rightarrow A_n \rightarrow 0.	\eqno(6)$$

\begin{lemma} \label{s7}
 For every $n\geq 1$ the restriction maps
$\rho_n\colon \mbox{Pic}(D(n)) \rightarrow \mbox{Pic}(D)$ and
 $\rho\colon \mbox{Pic}(\widehat D)
 \rightarrow \mbox{Pic}(D)$ are bijective.
\end{lemma}

\noindent{\em Proof}.
 Since $h^2(D,J^n/J^{n+1}) = 0$ for every $n\geq 1,$
 to obtain the surjectivity of $\rho_n$ it is sufficient to copy $\cite{H},$
 part 2 of the proof of theorem 3.1 on  page 179.
The last assertion follows from the bijectivity of all maps 
$\rho_n,$ just by the definition of line bundle on a formal
 scheme.\hfill\square{5}

\vspace{5mm}

\noindent The proof of Lemma \ref{s7} gives without any modification the 
analogous result with $S$ in place of $D.$

\begin{lemma} \label{s8}
 For every  $n\geq 1$ the restriction maps 
$\rho_n\colon \mbox{Pic}(S(n)) \rightarrow \mbox{Pic}(S)$ and 
$\rho\colon \mbox{Pic}(\widehat{S}) 
\rightarrow \mbox{Pic}(S)$ are bijective.
\end{lemma}

\noindent By Lemma \ref{s7} we can write 
${\cal O}_{D(n)}(a)$ (resp. ${\cal O}_{\widehat{D}}(a))$ 
for the unique line bundle on $D(n)$ (resp. $\widehat{D})$
 whose restriction to $D$ has degree $a.$ 
By Lemma \ref{s8} we can write ${\cal O}_{S(n)}(a,b)$ 
(resp. ${\cal O}_{\widehat{S}}(a,b))$ for the unique 
(up to  isomorphism) line bundle on $S(n)$ (resp. $\widehat{S}$
 whose  restriction to S is isomorphic to ${\cal O}_S(a,b).$

\begin{guess} \label{s9}
If the restriction 
of $G_1$ to a fiber $D$  
has splitting type $(a,\cdots ,a,a-1,\cdots ,a-1),$
 with $a-1$ appearing $x$ times,
 then $$G|_{\widehat{D}} \simeq 
{\cal O}_{\widehat{D}}(a)^{\oplus (r-x)}\oplus {\cal O}_{\widehat{D}}
(a-1)^{\oplus x }.$$
\end{guess}

\noindent{\em Proof}.
 Let $J$ be the ideal sheaf of $D$ in $V$ (or $ \widehat{S}).$
 We have $J/J^2 \simeq {\cal O}_D(t)\oplus{\cal O}_D.$
Since by our assumptions $t>0,$ then for every integer $n\geq 1$ the 
sheaf of  $ {\cal O}_D-$modules  $J^n/J^{n+1} \simeq {S^n}(J/J^2)$
 is the direct sum of $n+1$ line bundles on $D$ with nonnegative degree.
 Set $$A_n\colon= Hom({\cal O}_{D(n)}(a)^{\oplus(r-x)}
\oplus {\cal O}_{D(n)}(a-1)^{\oplus x},G_n).$$
 Thus $\{A_n\}$  is a rank $r^2$ vector bundle on
$\widehat{S}$ and $h^1(D,A_1) = 0.$
 Fix an integer $n\geq  1$ and assume $G_n
 \simeq {\cal O}_{D(n)}(a)^{\oplus(r-x)}
\oplus {\cal O}_{D(n)}(a-1)^{\oplus x}.$
 Fix $m_n \in 
H^0(D(n),A_n)$ with $m_n$ invertible. 
We have $H^1(D,J^n/J^{n+1}\otimes (a-1)) = 0.$
 Thus by $(5)$ we may lift $m_n$ to $m_{n+1} \in
H^0(D(n+1),A_{n+1})$ with $m_{n+1}|_{D(n)} = m_n.$
 By Nakayama's lemma $m_n$ is invertible.\hfill \square{5}

\subsection{Case 2: the rank divides the degree}
\label{10}

Fix integers $r$ and $a$ with $r\geq 2 .$
 Let $F$ be a simple rank $r$ vector bundle on $S$ 
such that $\mbox{det}(F)\cdot {\cal O}_S({\bf h}) = ar,$ i.e.,
 such that the restriction of to any fiber of $u$ has degree $d=ar.$
 Since $F$ is simple, by \ref{4.1} we have that
 $F$ is a flat limit of a flat 
family of vector bundles on $S.$ The general element of this 
family is a
 vector bundle $G$ whose restriction to a general fiber
 of $u$ has splitting type $(a,\cdots ,a)$ 
(i.e. it is isomorphic to 
${\cal O}_{{\bf_P}^1}(a)^{\oplus r})$
 but for which there are finitely many (say $z)$
 fibers of $u$ such that the restriction of $G$
 to each of these $z$ fibers has splitting type 
$(a+1,a,\cdots ,a,a-1)$ 
(i.e. it is isomorphic to 
${\cal O}_{{\bf_P}^1}(a+1)
\oplus{\cal O}_{{\bf_P}^1}(a)^{\oplus(r-2)}
\oplus{\cal O}_{{\bf_P}^1}(a-1).$
 The $z$ fibers of $u$ arising in this way are
 called the { \it jumping fibers} of $G.$

\vspace{5mm}
\noindent{\bf Theorem C} {\it  
Let $z$ be number of  jumping fibers of $G.$
Set  $E=G(-a{\bf h})$ and $m\colon=  \mbox{deg}(u_*E).$ Then 
$$z=c_2 (E)=c_2(G)-a(r-1)c_1(G)\cdot {\bf h}-ea^2r(r-1)/2$$
 and 
$$m=c_1(u_*E)=-z+c_1(G)\cdot{\bf h}+rae.$$}

\vspace{5mm}
\noindent{\em Proof}. If $E=G(-a{\bf h})$
then the restriction of $E$ to a general fiber has splitting type
$(0,\cdots, 0)$ whereas  the jumping fibers of $E$ have type
$(1,0, \cdots, 0 , -1).$
It follows that 
for each point  $p \in C$ such that $u^{-1}(p)$ is a jumping fiber
of $E$  the length $l(R^1u_*E_p(-1))=1.$ Hence each jumping fiber
has multiplicity one and   
$z=\sum_p l(R^1u_*E_p(-1))=c_2(E).$
By\cite{H2}  Lemma 2.1 
$c_2(E)=c_2(G)-a(r-1)c_1(G)\cdot {\bf h}-ea^2r(r-1)/2.$

By Grothendieck--Riemann--Roch
$\mbox{ch}(u_!E)\mbox{td}(T_C)=u_*\left(\mbox{ch}(E)\mbox{td}(T_S)\right),$
which gives 
$$c_1(u_*E)=-z+u_*\left((c_1-ra{\bf h})\cdot({\bf h}+{e+2 \over2}
{\bf f})\right)$$
and it follows that 
$m=c_1(u_*E)=-z+c_1(G)\cdot{\bf h}+rae.$
\hfill\square{5}

\vspace{5mm}
 
Changing bases, the coherent sheaf 
$H\colon= u_*(G\otimes{\cal O}_S(-a{\bf h}))$
 is a rank $r$ vector bundle on ${\mathbb P}^1.$
 The sheaf $u^*(H)$ is a rank $r$ subsheaf of 
$G$ with $G/(u^*(H)\otimes{\cal O}_S(a{\bf h}))$
 supported at the jumping fibers of $G.$
 The integers $e, r, a$ and the Chern classes of $G$ are
 uniquely determine the integers $z$ and  deg$(H).$ 
Conversely, the Chern classes of $G$ are uniquely 
determined by the integers $e, r, a,$ deg$(H)$ and $z.$
$G$ is a flat limit of a family $\{G_{\alpha}\}$ 
of simple vector bundles with the same properties 
of splitting type with respect to the fibers of $u,$ 
but with the added condition that $H_{\alpha} =
 u^*(G_{\alpha}(a{\bf h}))$ is rigid.
 Hence for such $G_{\alpha}$ the vector bundle $H_{\alpha}$
 is uniquely determined by the integers $r$ and deg$(H).$ 
Only the position of the $z$ points of $G_{\alpha}/u^*(H_{\alpha})(a{\bf h})$
 and the extension class of the sheaf $G_{\alpha}/u^*(H_{\alpha})(a{\bf h})$
 by the sheaf $u^*(H_{\alpha})(a{\bf h})$  depend on ${\alpha}.$

\begin{guess} \label{11} Let $\{G_n\}$  be a rank $r$ vector bundle on
$\widehat{S}$ such that $G_1$ is simple. Let $D$ 
be a fiber of $u$ such that $G_1|_D$
has splitting type $(a,\cdots,a).$
Then $G|_{\widehat{D}} \simeq {\cal O}_{\widehat{D}}(a)^{\oplus r}.$
\end{guess}

\noindent{\em Proof}.
 Twisting by the line bundle ${\cal O}_{\widehat{D}}(a)$
 we reduce to the case $a = 0.$ Let $J$ be the ideal sheaf of $ D$ 
in  $\widehat{S}.$ We have $J/J^2 \simeq 
{\cal O}_D(t)\oplus{\cal O}_D$ 
and  for every integer $n\geq 1$ the ${\cal O}_D-$sheaf 
$J^n/J^{n+1} \simeq S^n(J/J^2)$
 is the direct sum of $n+1$ line bundles on $D$ 
with nonnegative degree. Set 
$$A_n\colon= Hom({\cal O}_D(n)^{\oplus r},G_n).$$
 Thus $ \{A_n\}$ is a rank $r^2$ vector bundle 
on $\widehat{S}$ and $A_1$ is trivial.
 Fix an integer $n\geq  1$ and assume $G_n$ is  trivial.
 Fix a trivialization of $G_n,$ 
i.e. fix $m_n \in H^0(D(n),A_n)$ with $m_n$ invertible. 
Then $H^1(D,J^n/J^{n+1}\otimes A_1) = 0$
and by $(5)$ we may lift $m_n$ to $m_{n+1}\in H^0(D(n+1),A_{n+1})$
 with $m_{n+1}|_{D(n)} = m_n.$ 
By Nakayama's lemma $m_n$ is invertible. \hfill\square{5}

\vspace{5mm}

\noindent{\em Proof of  Theorem A}: \label{s0} For smoothness 
of the deformation space apply Lemma \ref{14}
 and Remark 
\ref{2}. Properties $(\iota)$ and $(\iota\iota)$ follow from Propositions
\ref{s9} and \ref{11}. \hfill\square{5}

\section{Bundles on a neighborhood of a 
 surface ruled  over a curve of higher genus}
We re-study the 
theory just done, now for  the case in which 
the divisor $S$ is ruled over a smooth curve 
$C$ of genus $q > 0.$
$S\colon= {\mathbb P}(B)$ is the 
projectivization of  a rank two vector bundle $B$ over $C.$ 
Let $u\colon S \rightarrow C$ be the projection.
 Fix a section, $h,$ of $S$ with minimal self-intersection 
and set $e\colon = - h^2.$ Denote by ${\bf h}$ the class of $h.$
By a theorem of M. Nagata and C. Segre 
 $e\geq - q$ and every integer $\geq - q$ 
 occurs for some rank two vector bundle $B$ on $C$ 
(see the introduction of \cite{LN}).  Pic$(S) \simeq 
u^*(\mbox{Pic}(C))\oplus {\mathbb Z}[\bf h]$ and 
 $u^*(M)\cdot {\bf h} =
 \mbox{deg}(M)$ and $u^*(M)\cdot u^*(H) = 0$ for all 
$M,H \in \mbox{Pic}(C).$
 The ruling $u$ induces an isomorphism between  $h$ and 
$C$ and we use this isomorphism to identify the normal bundle
 of $h$ in $S$ with a line bundle $ N$ on $C$ 
with deg$(N) = -e.$ The canonical line bundle of $S$
 is isomorphic to $u^*(w_C \otimes N)(-2 {\bf  h}).$
 We assume that $S$ is contained in a smooth threefold $W.$ 
We use the letter  $V$ to denote either a small  neighborhood of $S$
in the smooth topology  in $W,$ or else  the germ of $W$ around $V.$ 
Let $\widehat{S}$ be the formal completion of $S$ in $V$ and $I$ 
be the ideal sheaf of $S$ in  $V$ or in $\widehat{S}.$ 
Define ${\cal A} \in \mbox{Pic}(C)$ and $t \in {\mathbb Z}$ by the
 relation $I/I^2 \simeq u^*({\cal A})(t{\bf h}).$ We assume
 $t > 0$ and $deg({\cal A}) > 2q- 2 + |e|.$ By Riemann  Roch,
  $h^0(C,{\cal A}) > 0$ and $h^0(C,{\cal A} \otimes N) > 0.$

\begin{definition}{\em \label{def}
A  {\it good polarization} of a ruled surface  $u\colon S \rightarrow C$
is an ample divisor $R\in \mbox{Pic}(S)$ such that 
$R\cdot w_S+R \cdot D < 0$ for every fiber $D$ of $u.$
 } \end{definition}

\begin{remark}{\em \label{q0}
Good polarizations exist and occur quite frequently. For instance,
 choose any ample divisor $H \in \mbox{Pic}(S).$
 Then, for any fiber $ D$ of $u$ we have $D \cdot w_S = -2,$ 
and  $H(tD) \cdot w_S+H(tD) \cdot 
D < 0$  for every integer $t  >> 0.$ 
Furthermore  $H(tD)$ is ample for every $t\geq  0$ because 
$D$ is numerically effective. Thus, for $t >> 0$
 the line bundle $H(tD)$ is a  good polarization.
}\end{remark}

Our definition of good polarization is exactly 
the definition for which Lemma \ref{4.0}
and Proposition \ref{4.1} work for $q>0$ and for any 
$R$--stable vector bundle on $S.$ 
In what follows we assume that $R$ is a good polarization of $S.$

\begin{lemma} \label{q1} 
 Let $E =\{E_n\}$ 
 be a vector bundle on $\widehat{S}$ such that $E\vert_S$ is $R$--stable 
in the sense of Mumford and Takemoto. Then,
 $h^2(S(n),End(E_n)) = 0$ for all $n \geq 1.$
\end{lemma}

\noindent{\em Proof}. First assume $n = 1.$ 
Since $E\vert_S$ is $R$--stable and $ R\cdot w_S < 0,$ 
we have $h^0(S, Hom(E\vert_S,E\vert_S \otimes w_S)) = 0.$
 By Serre duality,  $$h^0(S, Hom(E\vert_S,E\vert_S \otimes w_S))
 = h^2(S,End(E\vert_S)),$$ concluding the case $n = 1.$ 
Now assume $n\geq 2$ and that the result is true for the integer
 $n-1,$ i.e., assume $h^2(S(n-1),End(E_{n-1})) = 0.$ 
Define  ${\cal A}$ by   $I/I^2 \simeq u^*({\cal A})(t{\bf h})$
as in the introduction of this section and set ${\cal N}\colon = I/I^2.$
Since dim$(S(n)) = 2$ we have $h^3(S(n),A) = 0$ 
for every coherent analytic sheaf $A$ on $S(n).$ 
Using $(2)$ for the integer $n-1$ and the vector
 bundle $End(E_n)$ instead of $E_n, $ together with the
inductive assumption, we see that $h^2(S(n),End(E_n)) = 0$ 
if $h^2(S,End(E\vert_S)\otimes {\cal N}^{\otimes n}) = 0.$
 Since $h^0(C,{\cal A}^{\otimes n}) > 0$ and $R$ is a good polarization,
 we have 
$R \cdot w_S \otimes {\cal N}^{\otimes n} < 0.$
 By Serre duality,
$$h^2(S,End(E\vert_S)\otimes {\cal N}^{\otimes n})
 = h^0(S,Hom(E\vert_S,E\vert_S \otimes w_S \otimes  ({\cal N}^*)^{\otimes n})) $$
 which vanishes, by the $R$--stability of $E\vert_S.$ \hfill\square{5}

\vspace{5mm}
The proof of Lemma \ref{4.0} works verbatim in 
this situation just assuming that 
$F$ is $R$--stable for some good polarization $R$ on $S.$
Lemmas
\ref{s7} and \ref{s8} work 
verbatim.
Therefore, we can write ${\cal O}_{\widehat{D}}(a))$ 
for the unique line bundle on  $\widehat{D}$ 
whose restriction to $D$ has degree $a$ and 
we  write 
$ {\cal O}_{\widehat{S}}u^*(M)(t{\bf h}) $
for the unique (up to  isomorphism) line bundle on $\widehat{S}$ 
whose restriction to $S$ is isomorphic to ${\cal O}_S u^*(M)(t{\bf h}).$
We may now re-state 2.9 in the situation.

\begin{guess}  \label{q9} Fix integers 
$r, \,a$ and $x$ with $r\geq 2$ and $0 < x < r.$ 
Let $G = \{G_n\}$  be a vector bundle on
$\widehat{S}$ such that the restriction of $G_1$ to a fiber $ D$ 
of $u$ has  splitting type $(a,\cdots ,a,a-1,\cdots ,a-1),$ 
with $a-1$ appearing $x$ times. 
Then $G|_{\widehat{D}} \simeq
 {\cal O}_{\widehat{D}}(a)^{\oplus(r-x)}
\oplus {\cal O}_{\widehat{D}}(a-1)^{\oplus x }.$\end{guess}

Let $F$ be an $R$--stable rank $r$ 
vector bundle on $S$ such that $\mbox{det}(F)\cdot {\cal O}_S({\bf h}) = ar,$
that is, such that the restriction of to any 
fiber of $u$ has degree $ar.$ Then
 \ref{1} and \ref{3} hold
 with $E\vert_S$ $R$--stable in place of $E\vert_S$ simple, we obtain.

\begin{guess} \label{q11}
Fix integers $r,\, a$ with $r\geq a.$ Let $\{G_n\}$
  be a rank $r$ vector bundle on $\widehat{S}$ such that $G_1$ 
is $R$--stable. Let $D$ be a fiber of $u$ such that $G_1|_D
 \simeq {\cal O}_D(a)^{\oplus r}.$ 
Then $G|_{\widehat{D}}
\simeq {\cal O}_{\widehat{D}}(a)^{\oplus r}.$ 
\end{guess}

\noindent{\em Proof of Theorem B}: For smoothness
of the deformation space apply Lemma \ref{q1}
 and Remark 
\ref{2}. Properties $(\iota)$ and $(\iota\iota)$ follow from Propositions
\ref{q9} and \ref{q11}. \hfill\square{5}

{\small
\noindent Edoardo Ballico,
Dept. of Mathematics, University of Trento,
38050 Povo, Italy\\
 ballico@science.unitn.it, fax: 39-(046) 1881624 \\

\noindent Elizabeth Gasparim,
Dept. of Mathematics, New Mexico State University\\
Las Cruces NM, 88003, USA,
gasparim@nmsu.edu\\

\end{document}